\documentclass{ifacconf}
\usepackage{natbib,amsmath,amssymb,latexsym,epsfig,graphicx}
\renewenvironment{pf}{\smallbreak\noindent{\it Proof. }}{\hfill$\Box$\smallbreak}
\renewenvironment{pf*}[1]{\smallbreak\noindent{\it #1}}{\hfill$\Box$\smallbreak}
\newcounter{definition}

\newcounter{remark}
\newenvironment{rmk}{\addtocounter{remark}{1}\smallbreak\noindent
  {\em Remark \theremark.}}{\smallbreak}
\newcounter{example}
\newenvironment{ex}{\addtocounter{example}{1}\smallbreak\noindent
  {\bf Example \theexample.}}{\hfill$\Box$\smallbreak}
\newenvironment{ex*}[1]{\addtocounter{example}{1}\smallbreak\noindent
  {\bf Example \theexample{} --- {\bf #1}}}{\hfill$\Box$\smallbreak}
\newcommand{\realR}{\mathbb{R}}

\newcommand{\itemi}{($i$)}
\newcommand{\ii}{($ii$)}
\newcommand{\iii}{($iii$)}
\newcommand{\iv}{($iv$)}

\DeclareMathOperator{\diag}{diag}

\DeclareMathOperator{\sat}{sat}

\begin{document}

\begin{frontmatter}
\title{\text{\!\!\!\!\!\!\!\!\!\!\!\!\!\!\!\!\!\!\!\!\!\!\!\!\!\!\!\!\!\! Minimax Adaptive Control for State Matrix with Unknown Sign}} 

\thanks[footnoteinfo]{The author is a member of the excellence center ELLIIT. Financial support was obtained from the Swedish Research Council and the European Research Council (ERC) Advanced Grant No.834142 (ScalableControl). The work was also partially supported by the Wallenberg AI, Autonomous Systems and Software Program (WASP) funded by the Knut and Alice Wallenberg Foundation.}

\author{Anders Rantzer} 
\address{Lund University, Sweden (e-mail: rantzer@control.lth.se).}

\begin{abstract}
For linear time-invariant systems having a state matrix with uncertain sign, we formulate and solve a minimax adaptive control problem as a zero sum dynamic game. Explicit expressions for the optimal value function and the optimal control law are given in terms of a Riccati equation. The optimal control law is adaptive in the sense that past data is used to estimate the uncertain sign for prediction of future dynamics. Once the sign has been estimated, the controller behaves like standard $H_\infty$ optimal state feedback.
\end{abstract}
\end{frontmatter}

\vglue 10mm
\section{Introduction}
The history of adaptive control dates back at least to aircraft autopilot development in the 1950s. Following the landmark paper \cite{ast+wit73}, a surge of research activity during the 1970s derived conditions for convergence, stability, robustness and performance under various assumptions. For example, \cite{ljung77ana} analysed adaptive algorithms using averaging, \cite{goodwin1981discrete} derived an algorithm that gives mean square stability with probability one, while \cite{guo1995convergence} analyzed the optimal asymptotic rate of convergence. On the other hand, conditions that may cause instability were studied in \cite{ega79book}, \cite{ioannou1984instability} and \cite{rohrs1985robustness}. Altogether, the subject has a rich history documented in numerous textbooks, such as \cite{aastrom2013adaptive}, \cite{goodwin2014adaptive}, 
\cite{sastry2011adaptive} and \cite{astolfi2007nonlinear}. In this paper, the focus is on worst-case models for disturbances and uncertain parameters, as discussed in \cite{cusumano1988nonlinear,sun1987theory,meg+03}. The ``minimax adaptive'' paradigm was introduced for linear systems in \cite{didinsky1994minimax} and nonlinear systems in \cite{pan1998adaptive}. 

This is a update to an earlier version published as \cite{RantzerIFAC20}. That version had an error in the main result, which has now been corrected.
The outline of the paper is still the same: Sections~2 introduces notation. Section~3 states the problem and reformulates it as a zero sum dynamic game on standard form. The main results are presented in section~4 together with an example. Concluding remarks are given in section~5, while some technical proofs are deferred to an appendix.

\newpage
\vglue 10mm
\section{Notation}
The set of $n\times m$ matrices with real coefficients is denoted $\realR^{n\times m}$. The transpose of a matrix $A$ is denoted $A^\top$. For a symmetric matrix $A\in\realR^{n\times n}$, we write $A\succ0$ to say that $A$ is positive definite, while $A\succeq0$ means positive semi-definite. For $A, B\in\realR^{n\times m}$, the expression $\langle A,B\rangle$ denotes the trace of $A^\top B$. Given $x\in\realR^n$ and $A\in\realR^{n\times n}$, the notation $|x|^2_A$ means 
$x^\top Ax$. Similarly, given $B\in\realR^{m\times n}$ and $A\in\realR^{n\times n}$, the trace of $B^\top AB$ is denoted $\|B\|^2_A$.
For $y\in\realR$, define $\sat(y)$ to be $1$ if $y>1$ and $-1$ if $y<-1$, otherwise equal to $y$. 

\goodbreak

\section{Minimax Adaptive Control}
\label{sec:minmax}

This paper is devoted to the following problem: 

Let $Q\in\realR^{n\times n}$ and $R\in\realR^{m\times m}$ be positive definite matrices and let $B\in\realR^{n\times m}$. Given $A\in\realR^{n\times n}$,
and a number $\gamma>0$, find, if possible, a control law $\mu$ that for every initial state $x_0$ attains the infimum 
\begin{align}
  \inf_\mu\sup_{w,i,N}\sum_{t=0}^N\left(|x_t|_Q^2+|u_t|_R^2-\gamma^2|w_t|^2\right),
\label{eqn:infsup}
\end{align}
where $i\in\{-1,1\}$, $w_t\in\realR^n$, $N\ge0$ and the sequences $x$ and $u$ are generated according to
\begin{align}
  x_{t+1}&=iAx_t+Bu_t+w_t& t&\ge0 
\label{eqn:plant}\\
  u_t&=\mu_t(x_0,\ldots,x_t,u_0,\dots,u_{t-1}).
\label{eqn:mu_LQ}
\end{align}

The problem can be viewed as a dynamic game, where the $\mu$-player tries to minimize the cost, while the $(w,i)$-player tries to maximize it. If it wasn't for the parameter $i$, this would be the standard game formulation of $H_\infty$ optimal control \cite{Basar/B95}. In our formulation, the maximizing player can choose not only $w$, but also the parameter $i$. This parameter is unknown, but constant, so an optimal feedback law tends to ``learn'' the value of $i$ in the beginning, in order to exploit this knowledge later. 
Such nonlinear adaptive controllers can stabilize and optimize the behavior also when no linear controller can simultaneously stabilize (\ref{eqn:plant}) for both $i=1$ and $i=-1$. 

To accommodate the uncertainty in $i$ when deciding $u_t$, we will see that it is sufficient for the controller to consider historical data collected in the matrix
\begin{align}
  Z_t&=\sum_{\tau=0}^{t-1}{\begin{bmatrix}Bu_\tau-x_{\tau+1}\\x_\tau\end{bmatrix}
  \begin{bmatrix}Bu_\tau-x_{\tau+1}\\x_\tau\end{bmatrix}^\top}.
\label{eqn:Z_t}
\end{align}
This gives
$
  \big\|\begin{bmatrix}I\;\;\,iA\end{bmatrix}^\top\big\|^2_{Z_t}
  =\sum_{\tau=0}^{t-1}|w_\tau|^2
$ and our problem can be reformulated as follows:

Given $Q\succ0, R\succ0$, $\gamma>0$ and a system
\begin{align}
  \left\{\begin{array}{ll}
    x_{t+1}=v_t\\
    Z_{t+1}=Z_t+{\begin{bmatrix}Bu_t-v_t\\x_t\end{bmatrix}
             \begin{bmatrix}Bu_t-v_t\\x_t\end{bmatrix}^\top},\quad &Z_0=0,
  \end{array}\right.
\label{eqn:xZ}
\end{align}
find, if possible, a control law 
\begin{align}
  u_t&=\eta(x_t,Z_t)
\label{eqn:eta}
\end{align}
that attains the infimum
{\begin{align}
  \!\!\!\!\!\inf_\eta\sup_{v,N}\Bigg[\sum_{t=0}^N\left(|x_t|_Q^2+|u_t|_R^2\right)
  -\gamma^2\min_i\big\|\begin{bmatrix}I\;\;\,iA\end{bmatrix}^\top\big\|^2_{Z_{N+1}}\Bigg]
\label{eqn:infsup2}
\end{align}
}when $x,u,Z$ are generated from $v$ and $x_0$ using
 (\ref{eqn:xZ})-(\ref{eqn:eta}).
\smallskip

In this formulation, the unknown sign $i$ does not appear in the dynamics, only in the penalty of the final state. As a consequence, no past states are needed in the control law (\ref{eqn:eta}), only the state $(x_t,Z_t)$. In fact, the problem is a standard zero-sum dynamic game \cite{basar1999dynamic}, which will next be addressed by dynamic programming.

\section{Minimax Dynamic Programming}

Define the Bellman operator $V\mapsto\mathcal{F}V$ by
\begin{align*}
  &\mathcal{F}V(x,Z):=\\
     &\min_u\max_v\left\{|x|_Q^2+|u|_R^2
     +V\biggl(v,Z+{\begin{bmatrix}Bu-v\\x\end{bmatrix}
     \begin{bmatrix}Bu-v\\x\end{bmatrix}^\top}\biggr)\right\}
\end{align*}
Then the following results holds:

\begin{thm}
  Given $A,B,Q,R$, define the operator $\mathcal{F}$ as above and $V_0,V_1,V_2\ldots$ according to the iteration
  \begin{align}
    V_0(x,Z)&=-\gamma^2\min_{i=\pm1}\big\|\begin{bmatrix}I\;\;\,iA\end{bmatrix}^\top\big\|^2_Z\label{eqn:Bellman_init}\\
    V_{k+1}(x,Z)&={\mathcal{F}}V_k(x,Z)\label{eqn:Bellman_t}
  \end{align}
  for $x\in\realR^n$ and positive semi-definite $Z\in\realR^{n\times n}$.
  The expressions (\ref{eqn:infsup}) and (\ref{eqn:infsup2}) have finite values if and only if the sequence $\{V_k(x,0)\}_{k=0}^\infty$ is upper bounded, in which case the limit $V_*:=\lim_{k\to\infty}V_k$ exists and $V_*(x_0,0)$ is equal to the values of (\ref{eqn:infsup}) and (\ref{eqn:infsup2}). 
  Defining $\eta(x,Z)$ as the minimizing value of $u$ in the expression for $\mathcal{F}V_*(x,Z)$ gives an optimal $\eta$ for (\ref{eqn:infsup2}), while the control law $\mu$ defined by 
  \begin{align}
    &\mu_t(x_0,\ldots,x_t,u_0,\dots,u_{t-1})\notag\\
    &=\eta\left(x_t,\sum_{\tau=0}^{t-1}{\begin{bmatrix}Bu_\tau-x_{\tau+1}\\x_\tau\end{bmatrix}
    \begin{bmatrix}Bu_\tau-x_{\tau+1}\\x_\tau\end{bmatrix}^\top}\right)
  \label{eqn:muopt}
  \end{align}
  is optimal for (\ref{eqn:infsup}). 
\label{thm:intro}
\end{thm}

\begin{pf}
  First note that $V_1\ge V_0$, so the sequence $V_0,V_1,V_2,\ldots$ is monotonically non-decreasing.

  For any fixed $N\ge0$, the value of (\ref{eqn:infsup}) is bounded below by the expression
  \begin{align}
    \inf_\mu\sup_{w,i}\sum_{t=0}^N\left(|x_t|_Q^2+|u_t|_R^2-\gamma^2|w_t|^2\right),
  \label{eqn:infsup_N} 
  \end{align}
  where $i\in\{-1,1\}$, $w_t\in\realR^n$ and the sequences $x$ and $u$ are generated according to (\ref{eqn:plant})-(\ref{eqn:mu_LQ}). 
  The value of (\ref{eqn:infsup_N}) grows monotonically with $N$ and (\ref{eqn:infsup}) is obtained in the limit. A change of variables with $v_t:=x_{t+1}$ and $Z_t$ given by (\ref{eqn:Z_t}) shows that (\ref{eqn:infsup_N}) is equal to
  {\small\begin{align}
    \inf_\mu\sup_{v}\left[\sum_{t=0}^N\left(|x_t|_Q^2+|u_t|_R^2\right)-\gamma^2\min_i\big\|\begin{bmatrix}I\;\;\,iA\end{bmatrix}^\top\big\|^2_{Z_{N+1}}\right]
  \label{eqn:sup_vT}
  \end{align}
  }where $x,Z,u$ are generated by (\ref{eqn:xZ}) combined with (\ref{eqn:mu_LQ}). 
  Standard dynamic programming shows that the value of (\ref{eqn:sup_vT}) is $V_{N+1}(x_0,0)$, where $V_k$ is defined by (\ref{eqn:Bellman_init})-(\ref{eqn:Bellman_t}).
  This proves that (\ref{eqn:infsup}) has a finite value if and only if the sequence $\{V_k(x,0)\}_{k=0}^\infty$ is upper bounded and the limit $V_*(x_0,0)$ is equal to the value of (\ref{eqn:infsup}). It follows from the relationship $V_k(x,Z)\le V_k(x,0)$ for $Z\succeq0$, that the limit then exists for all $Z\succeq0$ and by Lemma~\ref{lem:Zdiag} in the Appendix the same holds for all symmetric $Z$.
    
  If (\ref{eqn:infsup2}) is finite, then (\ref{eqn:sup_vT}) is bounded above by (\ref{eqn:infsup2}), so also $V_*:=\lim_{k\to\infty}V_k$ is finite. Conversely, if $V_*$ is finite, we may define $\eta(x,Z)$ as a minimizing value of $u$ in the expression for $\mathcal{F}V_*(x,Z)$. Then define the sequence $W_0,W_1,W_2,\ldots$ recursively by $W_0:=V_0$ and 
  \begin{align*}
    W_{k+1}&(x,Z):=\max_v\biggl\{|x|^2_Q+|\eta(x,Z)|^2_R\\
    &+W_k\biggl(v,Z+{\begin{bmatrix}B\eta(x,Z)-v\\x\end{bmatrix}
    \begin{bmatrix}B\eta(x,Z)-v\\x\end{bmatrix}^\top}\biggr)\biggr\}.
  \end{align*}
  By dynamic programming,
  \begin{align*}
    &W_N(x,0)\\
    &=\sup_{v}\left[-\gamma^2\min_i\big\|\begin{bmatrix}I\;\;\,iA\end{bmatrix}^\top\big\|^2_{Z_{N+1}}\!+\!\sum_{t=0}^N\left(|x_t|_Q^2+|u_t|_R^2\right)\right],
  \end{align*}
  where $x,Z,u$ are generated by (\ref{eqn:xZ}) combined with (\ref{eqn:eta}). Hence (\ref{eqn:infsup2}) is bounded above by $\lim_{k\to\infty}W_k(x_0,0)$. The definitions of $V_k$ and $W_k$ give by induction $V_*\ge W_k\ge V_k$ for all $k$, so $\lim_{k\to\infty}W_k=V_*$. This proves that the value of (\ref{eqn:infsup2}) equals $V_*(x_0,0)$ and $\eta$ is a minimizing argument. 
\end{pf}

\medskip
\begin{cor}
  With notation as in Theorem~\ref{thm:intro}, suppose that
  (\ref{eqn:infsup}) has a finite value and let $V_*:=\lim_{k\to\infty}V_k$. Then the Riccati equation
  {\small\begin{align}
    |x|^2_{P}&=\min_u\max_v\left\{|x|_Q^2+|u|_R^2-\gamma^2|Ax+Bu-v|^2+|v|^2_{P}\right\}
  \label{eqn:mRic}
  \end{align}
  }has a solution $0\prec P\prec \gamma^2I$ and the sequence defined by 
  \begin{align}
    \bar{V}_0(x,Z)&=|x|^2_P-\gamma^2\min_{i=\pm1}\big\|\begin{bmatrix}I\;\;\,iA\end{bmatrix}^\top\big\|^2_Z\label{eqn:Bellman_bar}\\
    \bar{V}_{k+1}(x,Z)&={\mathcal{F}}\bar{V}_k(x,Z).\label{eqn:Bellman_bar_t}
  \end{align}
  satisfies 
  $\bar{V}_0\le \bar{V}_1\le \cdots\le \lim_{k\to\infty}\bar{V}_k =V_*$.
\label{cor:Vbar}
\end{cor}

\begin{pf}
    Suppose that (\ref{eqn:infsup}) has a finite value. By Theorem~\ref{thm:intro}, this implies that the sequence $\{V_k(x,0)\}_{k=0}^\infty$ defined by (\ref{eqn:Bellman_init})-(\ref{eqn:Bellman_t}) is upper bounded.
    Define
    \begin{align*}
      V_k^+(x,Z)&:=|x|^2_{P_k}-\gamma^2\big\|\begin{bmatrix}I\;\;\,A\end{bmatrix}^\top\big\|^2_Z\\
      V_k^-(x,Z)&:=|x|^2_{P_k}-\gamma^2\big\|\begin{bmatrix}I\;\;\,-A\end{bmatrix}^\top\big\|^2_Z,
    \end{align*} 
    where $P_0=0$ and $P_k$ is given by the Riccati recursion
    {\small\begin{align*}
      |x|^2_{P_{k+1}}&=\min_u\max_v\left\{|x|_Q^2+|u|_R^2-\gamma^2|Ax+Bu-v|^2+|v|^2_{P_k}\right\}.
    \end{align*}
    }Then $V_k(x,Z)\ge \max\{V_k^+(x,Z),V_k^-(x,Z)\}$ for all $k$. This is trivial for $k=0$ and follows by induction for $k>0$, since $\mathcal{F}V^+_{k+1}=V^+_k$ and $\mathcal{F}V^-_{k+1}=V^-_k$. In the limit, it follows that the limit $P=\lim_{k\to\infty}P_k$ exists and 
    \begin{align}
      V_*(x,Z)\ge \bar{V}_0(x,Z).
    \label{eqn:Pbound}
    \end{align}
    Repeated application of $\mathcal{F}$ gives $V_*=\lim_{k\to\infty}\bar{V}_k$.
\end{pf}

\section{An Explicit Optimal Control Law}

The following result specifies a minimax optimal adaptive controller on explicit form for a range of $\gamma$-values. 

\begin{thm}
  Given $A\in\realR^{n\times n}$, $B\in\realR^{n\times m}$ and some positive definite $Q\in\realR^{n\times n}$, $R\in\realR^{m\times m}$, 
  assume that (\ref{eqn:mRic}) has a solution $0\prec P\prec \gamma^2I$, with minimizing argument $u=-Kx$. Define the sequence $\{\bar{V}_k\}_{k=0}^\infty$ by (\ref{eqn:Bellman_bar})-(\ref{eqn:Bellman_bar_t}). Then the following two conditions are equivalent:
  \begin{itemize}
    \item[\itemi] The value of (\ref{eqn:infsup}) is finite and $\bar{V}_k=\bar{V}_1$ for $k\ge1$.\\
    \item[\ii] The matrix $T:=Q+A^\top(P^{-1}-\gamma^2I)^{-1}A$ satisfies
    \begin{align*}
      &(\gamma^2I-P)(T-P)^{-1}(\gamma^2I-P)\\
      &\succeq(I\pm BKA^{-1})(\gamma^2I-P)(I\pm BKA^{-1})^\top.
    \end{align*}
    for both the ``+'' and the ``-'' sign.
  \end{itemize}
  Moreover, if \itemi{} and \ii{} hold, 
  then (\ref{eqn:infsup}) has the minimizing control law
  \begin{align}
    u_t
    &=\sat\left(\frac{\gamma^2\sum_{\tau=0}^{t-1}(Bu_\tau-x_{\tau+1})^\top Ax_\tau}{|x_t|^2_{T-P}}\right)Kx_t
  \label{eqn:optimalcontrol}
  \end{align}
  and, with notation $Y:=\gamma^2\begin{bmatrix}I&0\end{bmatrix}Z\begin{bmatrix}0&I\end{bmatrix}^\top$, the limit in Theorem~\ref{thm:intro} can be written as 
    \begin{align*}
      &V_*(x,Z)\notag\\
      &=\begin{cases}
         |x|^2_{P}-\gamma^2\min_{i=\pm1}
               \left\|\begin{bmatrix}I&iA\end{bmatrix}\right\|_Z^2
               \qquad\quad\,\,\,\text{if }|\langle A,Y\rangle|\ge |x|^2_{T-P}\\
         |x|^2_{T}-\gamma^2\left\|\diag\{I,A\}\right\|_Z^2
               +\langle A,Y\rangle^2|x|^{-2}_{T-P}
               \quad\text{otherwise.}
       \end{cases}
    \end{align*}
\label{thm:main}
\end{thm}
 
\begin{rmk}
  The intuition behind the optimal control law in Theorem~\ref{thm:main} is simple: The cases $u_t=Kx_t$ and $u_t=-Kx_t$ describe the situation when historical data collected in the expression $\sum_{\tau=0}^{t-1}(Bu_\tau-x_{\tau+1})^\top Ax_\tau$ is rich enough to make a reliable estimate about the uncertain parameter $i$. This estimate is then used as truth and the corresponding $H_\infty$ state feedback control law is applied. In the intermediate case, the historical data does not give a conclusive answer, so the controller gain is down-scaled accordingly.
\end{rmk}



The proof will be based on Lemma~\ref{lem:Zdiag}-\ref{lem:CDM} in the Appendix.

\begin{pf}
First assume \ii. We will prove the inequality $\mathcal{F}\bar{V}_1(x,Z)\le\bar{V}_1(x,Z)$ for all $x,Z$, (which together with Corollary~\ref{cor:Vbar} implies $\mathcal{F}\bar{V}_k=\bar{V}_1$ also for $k\ge1$). By Lemma~\ref{lem:Zdiag}, we may restrict attention to $Z$ of the form
\begin{align*}
  Z&=\gamma^{-2}\begin{bmatrix}0&Y\\Y^\top&0\end{bmatrix}.
\end{align*}
Let $\hat{u}$ be defined as in Lemma~\ref{lem:Vbar1}. Then the desired inequality is proved by the following sequence of relationships:
{\begin{align*}
  &\mathcal{F}\bar{V}_1(x,Z)-|x|_Q^2\\
  &=\min_u\max_v\left\{|u|_R^2
       +\bar{V}_1\biggl(v,Z+{\begin{bmatrix}Bu-v\\x\end{bmatrix}
       \begin{bmatrix}Bu-v\\x\end{bmatrix}^\top}\biggr)\right\}\\
  &\le|\hat{u}|_R^2
       +\max_v\bar{V}_1\biggl(v,Z+{\begin{bmatrix}B\hat{u}-v\\x\end{bmatrix}
       \begin{bmatrix}B\hat{u}-v\\x\end{bmatrix}^\top}\biggr)\\
  &=\max_{v}\max_{|\theta|\le1}\Big\{|v|^2_T-\theta^2|v|^2_{T-P}-\gamma^2|B\hat{u}-v|^2-\gamma^2|Ax|^2\notag\\
  &\qquad\qquad\qquad\qquad\quad-2\gamma^2\theta(Ax)^\top(B\hat{u}-v)-2\theta\langle A,Y\rangle\Big\}\notag\\
  &=\max_{|\theta|\le1}\Big\{\gamma^4|\theta Ax+B\hat{u}|^2_{(\gamma^2I-T+\theta^2T-\theta^2P)^{-1}}\notag\\
  &\qquad\qquad-\gamma^2\big[|Ax|^2+|B\hat{u}|^2+2\theta(Ax)^\top B\hat{u}\big]-2\theta\langle A,Y\rangle\Big\}\notag\\
  &=\max_{i=\pm1}\Big\{\gamma^4|iAx+B\hat{u}|^2_{(\gamma^2I-P)^{-1}}\notag\\
  &\qquad\qquad-\gamma^2\big[|Ax|^2+|B\hat{u}|^2+2i(Ax)^\top B\hat{u}\big]-2i\langle A,Y\rangle\Big\}\notag\\
  &=\max_{i=\pm1}\Big\{|iAx+B\hat{u}|^2_{(P^{-1}-\gamma^{-2}I)^{-1}}
  -2i\langle A,Y\rangle\Big\}\\
  &=|\hat{u}|_R^2+\max_{v,i}\Big\{|v|^2_P+|iAx+B\hat{u}-v|^2
      -\gamma^2\big\|\begin{bmatrix}I\;\,iA\end{bmatrix}^\top\big\|^2_Z\Big\}\\
  &=|\hat{u}|_R^2+\max_v\bar{V}_0\left(v,Z+{\footnotesize\begin{bmatrix}B\hat{u}-v\\x\end{bmatrix}\begin{bmatrix}B\hat{u}-v\\x\end{bmatrix}^\top}\right)\\
  &=\bar{V}_1(x,Z)-|x|_Q^2.
\end{align*}
}The first equality is the definition of $\mathcal{F}$ and the following inequality is trivial. The second equality is an application of Lemma~\ref{lem:Vbar1}. The third equality follows by analytic elimination of $v$ (Schur complement). The fourth equality is the main one. It follows from Lemma~\ref{lem:CDM} with 
\begin{align*}
  C&:=(\gamma^{2}I-T)^{-1/2}A\\
  D&:=(\gamma^{2}I-T)^{-1/2}BK\\
  M&:=(\gamma^{2}I-T)^{-1/2}(T-P)(\gamma^{2}I-T)^{-1/2}\\
  c&:=-2\gamma^{-2}(Ax)^\top B\hat{u}-2\langle A,Y\rangle.
\end{align*}
Multiplying the matrices  
\begin{align}
  2I+M^{-1}+M-(I\pm DC^{-1})(I+M)(I\pm DC^{-1})^\top
\label{eqn:MDC}
\end{align}
by $(\gamma^{2}I-T)^{1/2}$ from left and right gives, after straightforward manipulations detailed in the Appendix, 
\begin{align}
  &(\gamma^2I-P)(T-P)^{-1}(\gamma^2I-P)\label{eqn:TPBKA}\\
  &-(I\pm BKA^{-1})(\gamma^2I-P)(I\pm BKA^{-1})^\top,\notag
\end{align}
so \ii{} in Lemma~\ref{lem:CDM} follows from the assumption that \ii{} holds in Theorem~\ref{thm:main}. The fifth and sixth equalities are trivial algebraic manipulations, while the seventh is the definition of $\bar{V}_0$. Finally, the eighth equality is the definition of $\bar{V}_1$, so the implication from \ii{} to \itemi{} in Theorem~\ref{thm:main} is proved. The expressions for the optimal control law and the optimal cost are given by Lemma~\ref{lem:AA}.

For the opposite implication, consider the same sequence of expressions as before, but with $\check{u}$ defined as the minimizing argument in the definition of $\mathcal{F}\bar{V}_1(x,Z)$. This gives
\begin{align*}
  &\mathcal{F}\bar{V}_1(x,Z)-|x|_Q^2\\
  &=\max_{|\theta|\le1}\Big\{\gamma^4|\theta Ax+B\check{u}|^2_{(\gamma^2I-T+\theta^2T-\theta^2P)^{-1}}\notag\\
  &\qquad\qquad-\gamma^2\big[|Ax|^2+|B\check{u}|^2+2\theta(Ax)^\top B\check{u}\big]-2\theta\langle A,Y\rangle\Big\}\notag\\
  &\ge\max_{i=\pm1}\Big\{\gamma^4|iAx+B\check{u}|^2_{(\gamma^2I-P)^{-1}}\notag\\
  &\qquad\qquad-\gamma^2\big[|Ax|^2+|B\check{u}|^2+2i(Ax)^\top B\check{u}\big]-2i\langle A,Y\rangle\Big\}\notag\\
  &=|\check{u}|_R^2+\max_v\bar{V}_0\left(v,Z+{\footnotesize\begin{bmatrix}B\check{u}-v\\x\end{bmatrix}\begin{bmatrix}B\check{u}-v\\x\end{bmatrix}^\top}\right)\\
  &\ge\bar{V}_1(x,Z)-|x|_Q^2.
\end{align*}
Assume that \ii{} fails. Then,
by Lemma~\ref{lem:CDM}, there exist $x$ and $Z$ that make the first inequality strict, so \itemi{} must fail too and the proof is complete.
\end{pf}

As a complement to the previous result, we also give a lower bound on the values of $\gamma$ for which a solution exists:

\begin{thm}
  With $A,B,P,Q,R,T$ as in Theorem~\ref{thm:main}, (\ref{eqn:infsup}) has no finite value unless $0\prec P\prec \gamma^2I$ and $T\preceq \gamma^2I$.
\label{thm:lower} 
\end{thm}

\begin{pf}
  Inserting the bound (\ref{eqn:Pbound}) into the right hand side of the Bellman equation gives\\[-1mm]
  {\small\begin{align*}
    &V_*(x,Z)-|x|_Q^2\\
    &\!\!\!\ge \mathcal{F}\bar{V}_0(x,Z)-|x|_Q^2\\
    &\!\!\!=\min_u\max_{v,i}\left\{|u|_R^2+|v|^2_P
    -\gamma^2|iAx+Bu-v|^2-\gamma^2\big\|\begin{bmatrix}I\;\;\,iA\end{bmatrix}^\top\big\|^2_Z\right\}\\
    &\!\!\!=\min_u\max_i\left\{|u|_R^2+|iAx+Bu|^2_S-\gamma^2\big\|\begin{bmatrix}I\;\;\,iA\end{bmatrix}^\top\big\|^2_Z\right\}\\
    &\!\!\!\ge|x|^2_{T-Q}-\gamma^2\big\|\diag\{I,A\}^\top\big\|^2_Z\\[-1mm]
  \end{align*}
  }where the second inequality follows from Lemma~\ref{lem:AA}.
  Inserting the new bound $V_*(x,Z)\ge|x|_T^2-\gamma^2\big\|\diag\{I,A\}^\top\big\|^2_Z$ into the Bellman equation in the same way gives
  \begin{align*}
    &V_*(x,0)\\
    &\ge\min_u\max_{v,i}\left\{|x|^2_Q+|u|_R^2+|v|^2_T
        -\gamma^2|Bu-v|^2-\gamma^2|Ax|^2\right\}
  \end{align*}
  The last inequality shows that $T\preceq\gamma^2I$, so the proof is complete.
\end{pf}

\begin{ex} 
  Consider now the case $n=m=Q=R=A=B=1$. First of all, Theorem~\ref{thm:lower} shows that the game has no finite value unless $\gamma\ge 2.01$. On the other hand,  
  Theorem~\ref{thm:main} gives an optimal strategy for the dynamic game (\ref{eqn:infsup}) whenever $\gamma\ge2.5232$. Specifically, consider the case  $\gamma=2.5232$. Solving the Riccati equation
gives $P=1.6985$, which is clearly in the interval $[0,\gamma^2]$.
It follows that $T=3.3165$ and condition \ii{} of Theorem~\ref{thm:main} marginally holds. For larger $\gamma$, the margin would be bigger.

An exact expression for the value function $V_*$ is now given by the formula in Theorem~\ref{thm:main}, which shows that
\begin{align*}
  &V_*\left(x,\begin{bmatrix}z_{11}&z_{12}\\z_{12}&z_{22}\end{bmatrix}\right)\\[3mm]
  &=\begin{cases}
    1.70x^2-6.37(z_{11}+z_{22}-2|z_{12}|)&\text{if }|z_{12}|\ge0.25x^2\\
    3.32x^2-6.37(z_{11}+z_{22})+\frac{25.05z_{12}^2}{x^{2}}&\text{otherwise}
  \end{cases}
\end{align*}
and the optimal control law is
  \begin{align*}
    u_t&=\sat\left(\frac{3.93\sum_{\tau=0}^{t-1}(u_\tau-x_{\tau+1})x_\tau}{x_t^2}\right)0.698x_t.
  \end{align*}
\end{ex}

\begin{figure}
\vskip -20mm 
\hskip 17mm\includegraphics[width=.65\hsize]{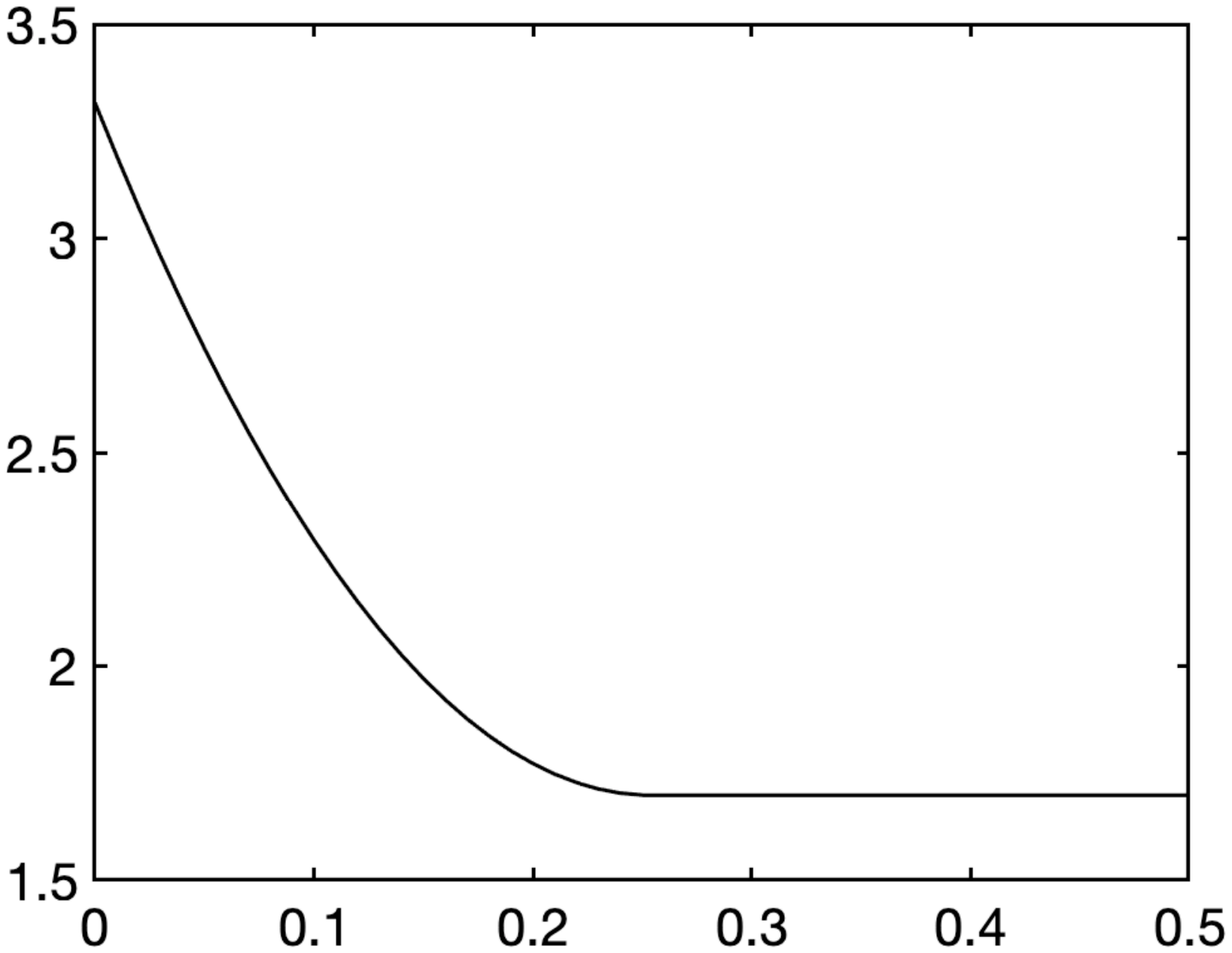}
\vskip -57mm $V_*\left(1,\begin{bmatrix}z&z\\z&z\end{bmatrix}\right)$
\vskip +25mm 
\hskip 7.2cm $z$
\caption{Optimal value function plot for Example~1. The variable $z$ represents information about the sign of $A$. The cost is maximal for $z=0$ (no information). For $|z|\ge0.25$, the cost is the same as if the sign was known.}
\label{fig:Vplot_a}  
\end{figure}

\goodbreak

\section{Concluding Remarks}
In this paper, we have formulated a control problem for uncertain linear systems as a zero-sum dynamic game. The solution is remarkable for two reasons:
\begin{enumerate}
  \item The dynamic programming formulation has an explicit solution in terms of a Riccati equation.\\
  \item The resulting optimal controller is adaptive: It reduces the aggressiveness of the controller until until enough data has been collected to get a parameter estimate that can be confidently trusted.
\end{enumerate}
The results are likely to be extendable to many other uncertainty structures. The case of uncertain input matrix $B$ will be particularly important, since the controller then needs to make active exploration in order to collect enough data for the exploitation phase.

\goodbreak

\bibliography{refs,rantzer}

\section*{Appendix: Supporting Lemmata}

The following lemmata are needed for the proof of Theorem~\ref{thm:main}.
\begin{lem}
  With the same notation as in Corollary~\ref{cor:Vbar}, it is true that 
  \begin{align*}
    \bar{V}_k(x,Z)&=\bar{V}_k(x,Z+\diag\{Z_{vv},Z_{xx}\})
  \end{align*}
  for every $x\in\realR^n$ and symmetric $Z\in\realR^{2n\times 2n}$, $Z_{vv}\in\realR^{n\times n}$ and $Z_{xx}\in\realR^{n\times n}$.
  \label{lem:Zdiag}
\end{lem}

\begin{pf}
  The proof is by induction. The statement is obviously true for $k=0$. Once it is proved for some value of $k$, the formula for $\mathcal{F}$ gives
  \begin{align*}
    &\bar{V}_{k+1}(x,Z+\diag\{Z_{vv},Z_{xx}\}-|x|_Q^2\\
    &=\min_u\max_v\biggl\{|u|_R^2
    +\bar{V}_k\biggl(v,Z+\diag\{Z_{vv},Z_{xx}\}\\
    &\qquad\qquad\qquad\qquad\qquad+{\begin{bmatrix}Bu-v\\x\end{bmatrix}
    \begin{bmatrix}Bu-v\\x\end{bmatrix}^\top}\biggr)\biggr\}\\
    &=\min_u\max_v\biggl\{|u|_R^2
    +\bar{V}_k\biggl(v,Z+{\begin{bmatrix}Bu-v\\x\end{bmatrix}
        \begin{bmatrix}Bu-v\\x\end{bmatrix}^\top}\biggr)\biggr\}\\
    &=\bar{V}_{k+1}(x,Z)-|x|_Q^2,
  \end{align*}
  so it holds also for $k+1$. Hence the result follows by induction over $k$.
\end{pf}

\begin{lem}
  Given matrices $A,B,Y$ and $P,Q,R,S\succ0$, suppose that
  \begin{align}
    |x|^2_P&=\min_u\left\{|x|_Q^2+|u|_R^2+|Ax+Bu|^2_S\right\},
    \label{eqn:QRS}
  \end{align}
  where the minimizing $u$ is given by $u=-Kx$ with
  \begin{align}
    K:=(R+B^\top SB)^{-1}B^\top SA.
  \label{eqn:K}
  \end{align}
  Then
  \begin{align*}
    &\min_u\max_{i\in\{-1,1\}}\left\{|x|_Q^2+|u|_R^2+|iAx+Bu|^2_S-2\langle iA,Y\rangle\right\}\notag\\[1mm]
    &=\begin{cases}
      |x|^2_P+2|\langle A,Y\rangle|&\text{if }|\langle A,Y\rangle| \ge|x|^2_{T-P}\\
      |x|^2_{T}+\langle A,Y\rangle^2|x|^{-2}_{T-P}&\text{otherwise}
    \end{cases} \\
    &=\max_{|\theta|\le1}\left\{|x|^2_T-\theta^2|x|^2_{T-P}-2\theta\langle A,Y\rangle\right\},
  \end{align*}
  where $T=Q+A^\top SA$. The maximizing $\theta$ is
  \begin{align*}
    \hat{\theta}&=-\sat\left(\frac{\langle A,Y\rangle}{|x|_{T-P}^2}\right)
  \end{align*}
  and the unique minimizing $u$ is $\hat{u}=-\hat{\theta}Kx$.
\label{lem:AA}
\end{lem}

\begin{pf}
  The definition of $K$ gives
  \begin{align}
    B^\top SA&=(R+B^\top SB)K.\notag
  \end{align}
  Multiplication by $K^\top$ from the left, and application of the identity
  \begin{align*}
    P=Q+K^\top RK+(A-BK)^\top S(A-BK)
  \end{align*}
  gives
  \begin{align*}
    K^\top B^\top SA=A^\top SBK=K^\top(R+B^\top SB)K={T-P}.
  \end{align*}
  The minimax theorem for convex-concave functions gives
  {\small\begin{align*}
    &\min_u\max_{i\in\{-1,1\}}\left\{|x|_Q^2+|u|_R^2+|iAx+Bu|^2_S-2\langle iA,Y\rangle\right\}\\  
    &\!\!\!=\max_{\theta_{-1},\theta_1}\min_u\sum_{i\in\{-1,1\}}\theta_i\left\{|x|_Q^2+|u|_R^2+|iAx+Bu|^2_S-2\langle iA,Y\rangle\right\}\\
    &\!\!\!=\max_{\theta\in[-1,1]}\min_u\left\{|x|_{Q+A^\top SA}^2+|u|_{R+B^\top SB}^2+2\theta(Ax)^\top SBu-2\theta \langle A,Y\rangle\right\}\\
    &\!\!\!=\max_{\theta\in[-1,1]}\left\{|x|^2_T-\theta^2|x|^2_{T-P}-2\theta \langle A,Y\rangle\right\}
  \end{align*}
  }where $\theta\in[-1,1]$, $\theta_{-1}=(1-\theta)/2$ and $\theta_{1}=(1+\theta)/2$. If $\langle A,Y\rangle \ge|x|^2_{T-P}$, the maximum over $\theta$ is attained by $\theta=-1$ and the value is $|x|^2_P+2\langle A,Y\rangle$.
  On the other hand, if $\langle A,Y\rangle \le -|x|^2_{T-P}$, the maximum is given by $\theta=1$ and the value is $|x|^2_P-2\langle A,Y\rangle$. Finally, if $|\langle A,Y\rangle|< |x|^2_{T-P}$, the optimal value of $\theta$ is in the interior of the interval $(-1,1)$ and determined by
  \begin{align*}
    |Ax-Bu|^2_S+2\langle A,Y\rangle
    &=|Ax+Bu|^2_S-2\langle A,Y\rangle\\
    \langle A,Y\rangle
    &=x^\top A^\top SBu\\
    &=-\theta x^\top A^\top SBKx\\
    &=-\theta|x|^2_{T-P}.
  \end{align*}
  This gives $u=|x|_{T-P}^{-2}\langle A,Y\rangle Kx$
  and the value
  \begin{align*}
    &\sum_i\theta_i\left\{|x|_Q^2+|u|_R^2+|iAx+Bu|^2_S-2\langle iA,Y\rangle\right\}\\
    &=|x|_Q^2+|u|_R^2+|Ax|^2_S+|Bu|^2_S\\
    &=|x|_Q^2+|Ax|^2_S+|\theta Kx|^2_{R+B^\top SB}\\
    &=|x|_{T}^2+\theta^2|x|^2_{T-P}\\
    &=|x|^2_{T}+\langle A,Y\rangle^2|x|^{-2}_{T-P}.
  \end{align*}
\end{pf}

\begin{lem}
  With the same notation as in Theorem~\ref{thm:main} and $Y:=\gamma^2\begin{bmatrix}I&0\end{bmatrix}Z\begin{bmatrix}0&I\end{bmatrix}^\top$, it is true that
  \begin{align*}
    &\bar{V}_1(x,Z)\notag\\
    &=\max_v\left\{|x|_Q^2+|\hat{u}|_R^2
    +\bar{V}_0\biggl(v,Z+{\begin{bmatrix}B\hat{u}-v\\x\end{bmatrix}
    \begin{bmatrix}B\hat{u}-v\\x\end{bmatrix}^\top}\biggr)\right\}\notag\\
    &=\max_{|\theta|\le1}\left\{|x|^2_T-\theta^2|x|^2_{T-P}-2\theta\langle A,Y\rangle-\gamma^2\big\|\diag\{I,A\}^\top\big\|^2_Z\right\}, 
  \end{align*}
  where $\hat{u}=-\hat{\theta}Kx$ and $\hat{\theta}$ is a maximizing argument of the last expression.
  \label{lem:Vbar1}
\end{lem}

\begin{pf}
  Define $S:=(P^{-1}-\gamma^{-2}I)^{-1}$. Then
  \begin{align*}
    &\bar{V}_1(x,Z)-|x|_Q^2\notag\\
    &=\inf_u\sup_v\left\{|u|_R^2
    +\bar{V}_0\left(v,{Z+\begin{bmatrix}Bu-v\\x\end{bmatrix}\!\!
      \begin{bmatrix}Bu-v\\x\end{bmatrix}^\top}\right)\right\}\notag\\
    &=\min_u\max_{v,i}\Big\{|u|_R^2+|v|_P^2-\gamma^2|iAx+Bu-v|^2\\
    &\qquad\qquad\qquad\qquad\qquad\qquad-\gamma^2\big\|\begin{bmatrix}I\;\;iA\end{bmatrix}^\top\big\|^2_Z\Big\}\notag\\
    &=\min_u\max_{i\in\{-1,1\}}\left\{|u|_R^2+|iAx+Bu|^2_S-\gamma^2\big\|\begin{bmatrix}I\;\;\,iA\end{bmatrix}^\top\big\|^2_Z\right\}.
  \end{align*}
  so the desired expressions follow from Lemma~\ref{lem:AA}.
\end{pf}

\begin{lem}
  Given $C,D,M\in\realR^{n\times n}$, with $C$ invertible and $M$ symmetric positive definite, the following statements are equivalent:
  \begin{itemize}
    \item[\itemi] For every $x\in\realR^n$ and $c\in\realR$ the maximum
    \begin{align*}
      \max_{\theta\in[-1,1]}\left(|\theta Cx+Dx|^2_{(I+\theta^2M)^{-1}}+\theta c\right)
    \end{align*}
    is attained for either $\theta=-1$ or $\theta=1$.
    \item[\ii] The two matrices
    \begin{align*}
      2I+M^{-1}+M-(I+DC^{-1})(I+M)(I+DC^{-1})^\top\\
      2I+M^{-1}+M-(I-DC^{-1})(I+M)(I-DC^{-1})^\top
    \end{align*}
    are positive semi-definite.
  \end{itemize}
\label{lem:CDM}
\end{lem}

\begin{pf}
To prove the lemma, it is convenient to consider two additional statements:
\begin{itemize}
  \item[\iii]For all $x\in\realR^n$ and $\theta\in[-1,1]$, it holds that
  \begin{align*}
    &|\theta Cx+Dx|^2_{(I+\theta^2M)^{-1}}\\
    &\le|Cx|^2_{(I+M)^{-1}}+|Dx|^2_{(I+M)^{-1}}
    +2\theta x^\top C^\top(I+M)^{-1}Dx.
  \end{align*}
  \item[\iv]For all $x\in\realR^n$ and $\theta\in[-1,1]$, it holds that
  \begin{align*}
      |\theta Cx+Dx|^2_{(I+M^{-1}+\theta^2I+\theta^2M)^{-1}}
      &\le|Cx|^2_{(I+M)^{-1}}.
  \end{align*}
\end{itemize}
We will first prove that \itemi{} is equivalent to \iii. Define
\begin{align*}
  f(\theta)&=|\theta Cx+Dx|^2_{(I+\theta^2M)^{-1}}\\
  g(\theta)&=|Cx|^2_{(I+M)^{-1}}+|Dx|^2_{(I+M)^{-1}}
      +2\theta x^\top C^\top(I+M)^{-1}Dx.
\end{align*}
Assuming that \iii{} holds gives
\begin{align*}
  \max_{\theta\in[-1,1]}[f(\theta)+\theta c]
  &\le \max_{\theta\in[-1,1]}[g(\theta)+c\theta]\\
  &=\max_{\theta=\pm1}[g(\theta)+c\theta]\\
  &=\max_{\theta=\pm1}[f(\theta)+\theta c],
\end{align*}
for all $x$, which proves \itemi. On the other hand, if \iii{} fails, there exist a $\hat{\theta}$ and $x$ such that $f(\hat{\theta})>g(\hat{\theta})$. With
\begin{align*}
  c&:=-2x^\top C^\top(I+M)^{-1}Dx
\end{align*}
it follows that
\begin{align*}
  f(\hat{\theta})+c\hat{\theta}
  &>g(\hat{\theta})+c\hat{\theta}\\
  &=|Cx|^2_{(I+M)^{-1}}+|Dx|^2_{(I+M)^{-1}}\\
  &=\max_{\theta=\pm1}[g(\theta)+c\theta]\\
  &=\max_{\theta=\pm1}[f(\theta)+\theta c],
\end{align*}
so also \itemi{} fails. Hence \itemi{} is equivalent to \iii.

Next, to prove that \iii{} is equivalent to \iv, note that 
\begin{align*}
  &(I+\theta^2M)^{-1}-(I+M)^{-1}\\
  &=(I+\theta^2M)^{-1}\left[I+M-(I+\theta^2M)\right](I+M)^{-1}\\
  &=(1-\theta^2)(I+\theta^2M)^{-1}M(I+M)^{-1}\\
  &=(1-\theta^2)(I+M^{-1}+\theta^2I+\theta^2M)^{-1}.
\end{align*}
Hence, by subtracting $|\theta Cx+Dx|^2_{(I+M)^{-1}}$ from both sides of the inequality in \iii{}, we get
\begin{align*}
  (1-\theta^2)|\theta Cx+Dx|^2_{(I+M^{-1}+\theta^2I+\theta^2M)^{-1}}
  &\le(1-\theta^2)|Cx|^2_{(I+M)^{-1}}.
\end{align*}
The equivalence between \iii{} and \iv{} follows.

It remains to prove that \iv{} is equivalent to \ii. The inequality in \iv{} holds for all $x$ if and only if
\begin{align*}
  &(\theta C+D)^\top(I+M^{-1}+\theta^2I+\theta^2M)^{-1}(\theta C+D)\\
  &\preceq C^\top(I+M)^{-1}C,
\end{align*}
or equivalently that the matrix
\begin{align*}
  I+M^{-1}+\theta^2I+\theta^2M-(\theta I+DC^{-1})(I+M)(\theta I+DC^{-1})^\top
\end{align*}
is positive semi-definite. The matrix is linear in $\theta$, so the condition holds for all $\theta\in[-1,1]$ if and only if it holds for $\theta=\pm1$. Equivalence between \iv{} and \ii{} follows.
\end{pf}

\begin{pf*}{Proof of the relationship between (\ref{eqn:MDC}) and (\ref{eqn:TPBKA}).}
  \begin{align*}
  &(\gamma^{2}I-T)^{1/2}(I+M)(\gamma^{2}I-T)^{1/2}\\
  &=\gamma^{2}I-T+T-P\\
  &=\gamma^2I-P\\\\
  &(\gamma^{2}I-T)^{1/2}(I+M^{-1})(\gamma^{2}I-T)^{1/2}\\
  &=\gamma^{2}I-T+(\gamma^{2}I-T)(T-P)^{-1}(\gamma^{2}I-T)\\
  &=(\gamma^{2}I-T)(T-P)^{-1}\big[(T-P)+(\gamma^{2}I-T)\big]\\
  &=(\gamma^{2}I-T)(T-P)^{-1}(\gamma^{2}I-P)\\\\
  &(\gamma^{2}I-T)^{1/2}(2I+M+M^{-1})(\gamma^{2}I-T)^{1/2}\\
  &=\big[(T-P)+(\gamma^{2}I-T)\big](T-P)^{-1}(\gamma^{2}I-P)\\
  &=(\gamma^{2}I-P)(T-P)^{-1}(\gamma^{2}I-P)\\\\
  &(\gamma^{2}I-T)^{1/2}(I\pm DC^{-1})(\gamma^{2}I-T)^{-1/2}=I\pm BKA^{-1}.
  \end{align*}
\end{pf*}
\end{document}